\documentclass% [12pt]{amsart}
[10pt]{amsart}
\usepackage{amssymb,epsfig,amscd}
\usepackage{xy}
\xyoption{all}
\usepackage[mathscr]{eucal}

\newtheorem{thm}{Theorem}

\newtheorem{conj}[thm]{Conjecture}
\numberwithin{equation}{subsection}
\numberwithin{thm}{subsection}

\newtheorem*{thm*}{Theorem}

\begin{document}

\newcommand\R{\mathbb{R}}
\newcommand\C{\mathbb{C}}
\newcommand\Z{\mathbb{Z}}
\newcommand\Q{\mathbb{Q}}
\newcommand\SL{\mathsf{SL}}
\newcommand\Sp{\mathsf{Sp}}
\newcommand\PSp{\mathsf{PSp}}
\newcommand\SLtR{\SL(2,\R)}
\newcommand\SLtC{\SL(2,\C)}
\newcommand\SLnR{\SL(n,\R)}
\newcommand\PSL{\mathsf{PSL}}
\newcommand\PSLtR{\PSL(2,\R)}
\newcommand\PSLtC{\PSL(2,\C)}
\newcommand\PGL{\mathsf{PGL}}
\newcommand\SO{\mathsf{SO}}
\newcommand\OO{\mathsf{O}}
\newcommand\SU{\mathsf{SU}}
\newcommand\Uo{\mathsf{U}(1)}
\newcommand\Un{\mathsf{U}(n)}
\newcommand\Uno{\mathsf{U}(n,1)}
\newcommand\Uoo{\mathsf{U}(1,1)}
\newcommand\SUoo{\mathsf{SU}(1,1)}
\newcommand\Upq{\mathsf{U}(p,q)}
\newcommand\Upp{\mathsf{U}(p,p)}
\newcommand\Uppq{\mathsf{U}(p,p)\times\mathsf{U}(q-p)}
\newcommand\SpnR{\Sp(n,\R)}
\newcommand\SpfR{\Sp(4,\R)}
\newcommand\tS{\tilde{S}}
\newcommand\trho{\tilde{\rho}}

\newcommand\rpt{\R\mathsf{P}^2}
\newcommand\rpo{\R\mathsf{P}^1}
\newcommand\rpn{\R\mathsf{P}^{n-1}}
\newcommand\rpthree{\R\mathsf{P}^3}
\newcommand\Ht{\mathsf{H}^2}
\newcommand\Hthree{\mathsf{H}^3}
\newcommand\SLthR{\mathsf{SL}(3,\R)}
\newcommand\PGLthR{\mathsf{PGL}(3,\R)}
\newcommand\cpo{\C\mathsf{P}^1}
\renewcommand\div{\mathsf{div}}
\newcommand\Isom{\mathsf{Isom^+}}
\newcommand{\Fricke}{\mathfrak{F}}
\newcommand{\FrickeX}{\Fricke(X)}
\newcommand{\CC}{\mathfrak{C}(\Sigma)}

\newcommand\hol{\mathsf{hol}}
\newcommand\Jac{\mathsf{Jac}(X)}
\newcommand\Mod{\mathsf{Mod}(\Sigma)}
\newcommand\tM{\tilde{M}}
\newcommand\End{\mathsf{End}}

\newcommand\Teich{\mathfrak{T}(\Sigma)}
\newcommand\TeichX{\mathfrak{T}(X)}
\newcommand\dev{\mathsf{dev}}
\renewcommand\deg{\mathsf{deg}}
\newcommand\rank{\mathsf{rank}}
\newcommand\GLnC{\mathsf{GL}(n,\C)}
\newcommand\HH{\mathcal{H}}

\newcommand\NS{\mathcal{M}_{n,0}(X)}
\newcommand\Aut{\mathsf{Aut}}
\newcommand\Inn{\mathsf{Inn}}
\newcommand\Out{\mathsf{Out}}
\newcommand\Rep{\mathsf{Rep}}
\newcommand\Diff{\mathsf{Diff}}

\newcommand\Euler{\mathsf{Euler}}
\newcommand\Hopf{\mathsf{Hopf}}

\newcommand\hK{hyper-K\"ahler\ } 

\newcommand\QF{\mathcal{QF}(\Sigma)}
\newcommand\Ep{E_{\tilde\rho}}

\newcommand\PU{\mathsf{PU}}
\newcommand\U{\mathsf{U}}

\newcommand\Hom{\mathsf{Hom}}
\newcommand\hp[1]{\Hom(\pi,#1)}
\newcommand\hpg{\hp{G}}
\newcommand\hpgg{\hpg/G}
\title[Higgs bundles and geometric structures]
{Higgs bundles and geometric structures on surfaces}

\author[W. Goldman]{William M.~Goldman}
\address{ Mathematics Department,
University of Maryland, College Park, MD  20742 USA  }
\email{ wmg@math.umd.edu }

\thanks{Goldman gratefully acknowledge partial support from 
National Science Foundation grants  
DMS-070781, DMS-0103889 and DMS-0405605. 
This paper was presented at the Geometry Conference honouring 
Nigel Hitchin, Consejo Superior de Investigaciones CientÃ­ficas,
Madrid, 5 September 2006. We gratefully acknowledge support from the Oswald Veblen
Fund at the Institute for Advanced Study where this work was completed.
}

\date{\today}
\dedicatory
{Dedicated to Nigel Hitchin on the occasion of his sixtieth birthday.}
\maketitle
\tableofcontents

\section*{Introduction}
In the late 1980's Hitchin~\cite{Hitchin_selfduality} and 
Simpson~\cite{Simpson1} 
discovered deep connections between representations of
fundamental groups of surfaces and algebraic geometry.
% (from a somewhat different point of view)   
The fundamental group $\pi = \pi_1(\Sigma)$ of a closed orientable 
surface $\Sigma$ of genus $g > 1$ is 
an algebraic object governing the topology of $\Sigma$.  For a Lie
group $G$, the space of conjugacy classes of representations $\pi\to G$
is a natural algebraic object $\hpgg$ whose geometry, topology and dynamics
intimately relates the topology of $\Sigma$ and the various geometries 
associated with $G$.
In particular $\hpgg$ arises as a moduli space of locally homogeneous
geometric structures as well as flat connections on bundles over
$\Sigma$.
 
Giving $\Sigma$ a conformal structure profoundly affects $\pi$ and its
representations.  This additional structure induces further geometric
and analytic structure on the deformation space $\hpgg$.  Furthermore
this analytic interpretation allows Morse-theoretic methods to compute
the algebraic topology of these non-linear finite-dimensional spaces.

For example, when $G = \Uo$, the space of representations 
is a torus of dimension $2g$. Give $\Sigma$ a conformal structure ---
denote the resulting Riemann surface by $X$.
The classical Abel-Jacobi theory 
identifies representations $\pi_1(X)\longrightarrow\Uo$ 
with topologically trivial holomorphic line bundles over $X$.
The resulting {\em Jacobi variety\/} is an abelian variety, whose structure
strongly depends on the Riemann surface $X$.
However the underlying 
{\em symplectic manifold \/} 
depends only on the topology of $\Sigma$, and indeed just
the fundamental group $\pi$.

Another important class of representations of $\pi$ 
arises from introducing the local structure of 
{\em hyperbolic geometry\/} to $\Sigma$. 
Giving $\Sigma$ a Riemannian metric of curvature $-1$ determines
a representation $\rho$ in the group $G = \Isom(\Ht) \cong \PSLtR$.
These representations, which we call {\em Fuchsian,\/}
are characterized as {\em embeddings\/}
of $\pi$ onto {\em discrete\/} subgroups of $G$.
Equivalence classes of Fuchsian representations comprise
the Fricke-Teichm\"uller space $\Fricke(\Sigma)$ 
of marked hyperbolic structures on $\Sigma$, which embeds
in $\hpgg$ as a connected component. This component is a cell 
of dimension $6g-6$ upon which
the mapping class group acts properly.

The theory of Higgs bundles, pioneered by Hitchin and Simpson,
provides an analytic approach to studying surface group
representations and their deformation space. The purpose of this paper
is to describe the basic examples of this theory, emphasizing 
relations to deformation and rigidity of geometric structures.
In particular we report on some very recent developments when $G$ is 
a real Lie group, either a split real semisimple group 
or an automorphism group of a Hermitian symmetric space of noncompact type.

In the twenty years since the appearance of Hitchin's and Simpson's work,
many other developments directly arose from this work. These relate  to
variations of Hodge structures, spectral curves, integrable systems, 
Higgs bundles over noncompact Riemann surfaces and higher-dimensional 
K\"ahler manifolds, and the finer topology of the deformation spaces. 
None of these topics are discussed here. It is an indication of the power and
the depth of these ideas that so mathematical subjects have been profoundly
influenced by the pioneering work of Hitchin and Simpson.

\subsubsection*{Acknowledgements}
I am grateful to Nigel Hitchin for the inspiration of these ideas, and
to him and Simon Donaldson for the opportunity to study with them at
the Maths Institute in Oxford in 1989.  Over the years my knowledge of
the subject has benefitted enormously from conversations with Yves
Benoist, Steve Bradlow, Marc Burger, Suhyoung Choi, Kevin Corlette,
Oscar Garcia-Prada, Peter Gothen, Olivier Guichard, Alessandra Iozzi,
Misha Kapovich, Fran\c cois Labourie, John Loftin, John Millson,
Ignasi Mundet i Riera, Walter Neumann, Carlos Simpson,
Ser-Peow Tan,
Richard Wentworth, Anna Wienhard, Graeme Wilkin, Mike Wolf, 
Scott Wolpert, and Eugene Xia. I also thank the students in the
Experimental Geometry Lab at the University of Maryland, in particular
Ryan Hoban, Rachel Kirsch and Anton Lukyanenko for help producing the
illustrations for this paper. I am grateful to the anonymous referee
for several helpful suggestions. Finally I wish to thank the Mathematical
Sciences Research Institute, the Institute for Advanced Study and the
Mathematics Department at Princeton University for their hospitality
where this work was completed.

\section{Representations of the fundamental group}

\subsection{Closed surface groups}
\begin{figure}[!ht]
\centerline{\epsfig{figure=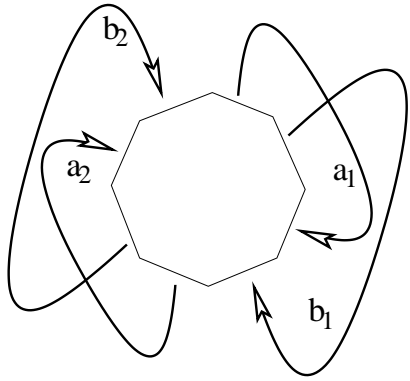,height=5cm}}
\caption{The pattern of identifications for a genus 2 surface.  The
sides of an octagon are pairwise identified to construct a surface of
genus 2. The 8 vertices identify to a single 0-cell in the quotient,
and the 8 sides identify to four 1-cells, which correspond to the four
generators in the standard presentation of the fundatmental group.}
\label{fig:octagon}
\end{figure}

\begin{figure}[ht]
\centerline{\epsfig{figure=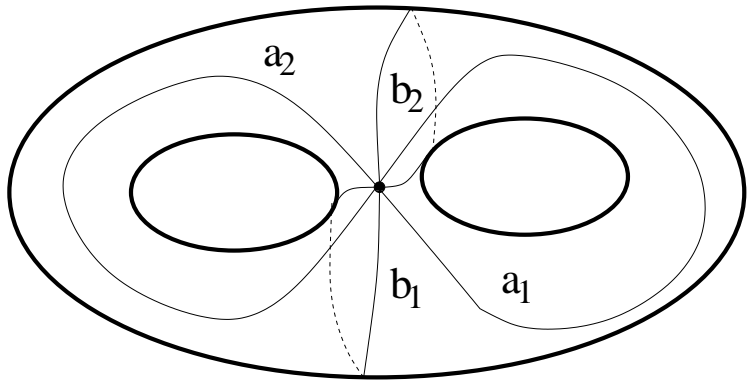,height=5cm}}
\caption{The genus 2 surface as an identification space.}
\label{fig:genustwo}
\end{figure}

Let $\Sigma = \Sigma_g$ be a closed orientable surface of genus
$g > 1$. Orient $\Sigma$, and choose a smooth
structure on $\Sigma$. Ignoring basepoints, denote the fundamental
group $\pi_1(\Sigma)$ of $\Sigma$ by $\pi$.
The familiar decomposition of $\Sigma$ as a 
$4g$-gon with $2g$ identifications 
(depicted in Figures~\ref{fig:octagon} and~\ref{fig:genustwo})
of its sides leads to a presentation
\begin{equation}\label{eq:presentation}
\pi  =   \langle A_1,B_1,\dots, A_g,B_g \mid 
[A_1,B_1]\dots[A_g,B_g]  = 1 \rangle 
\end{equation}
where $[A,B]:= ABA^{-1}B^{-1}$. 

% \clearpage
\subsection{The representation variety}

Denote the set of representations $\pi\xrightarrow{\rho} G$ by $\hpg$.
Evaluation on a collection $\gamma_1,\dots,\gamma_N\in\pi$ defines a map:
\begin{align} \label{eq:eval}
\hpg & \longrightarrow G^N \notag\\
\rho & \longmapsto \bmatrix \rho(\gamma_1)  \\ \vdots \\ \rho(\gamma_N)
\endbmatrix 
\end{align}
which is an embedding if $\gamma_1,\dots,\gamma_N$ generate $\pi$.  Its
image consists of $N$-tuples 
\begin{equation*}
(g_1,\dots,g_N)\in G^N  
\end{equation*}
satisfying equations
\begin{equation*}
R(g_1,\dots,g_N) = 1 \in G 
\end{equation*}
where $R(\gamma_1,\dots,\gamma_N)$ are defining relations in $\pi$
satisfied by $\gamma_1,\dots,\gamma_N$. 
If $G$ is
a linear algebraic group, 
these equations are polynomial equations in the matrix entries of $g_i$.
Thus the evaluation map \eqref{eq:eval}
identifies $\hpg$ as an algebraic subset of
$G^N$. The resulting algebraic structure is independent of the generating
set. In particular $\hpg$ inherits both the Zariski and the classical
topology. We consider the classical topology unless
otherwise noted.

In terms of the standard presentation \eqref{eq:presentation}, $\hpg$
identifies with the subset of $G^{2g}$ consisting of
\begin{equation*}
(\alpha_1,\beta_1,\dots \alpha_g,\beta_g)
\end{equation*}
satisfying the single $G$-valued equation
\begin{equation*}
[\alpha_1,\beta_1]\dots[\alpha_g,\beta_g] = 1.
\end{equation*}

\subsection{Symmetries}
The product $\Aut(\pi)\times\Aut(G)$ acts
naturally by left- and right-composition,
on $\hpg$:
An element 
\begin{equation*}
(\phi,\alpha)\in \Aut(\pi)\times\Aut(G)  
\end{equation*}
transforms $\rho\in\hpg$ to the composition
\begin{equation*}
% \alpha \circ \rho \circ \phi^{-1}:\;
\pi \xrightarrow{\phi^{-1}}
\pi \xrightarrow{\rho}
G \xrightarrow{\alpha} G.
\end{equation*}
The resulting action 
preserves the algebraic structure on $\hpg$ .

\subsection{The deformation space}

For any group $H$, let $\Inn(H)$ denote the normal subgroup of
$\Aut(H)$ comprising {\em inner autormorphisms.\/}  The quotient group
$\Aut(H)/\Inn(H)$ is the {\em outer automorphism group,\/} denoted
$\Out(H)$.

We will mainly be concerned with the quotient
\begin{equation*}
\hpgg := \hpg/ \big( \{1\}\times \Inn(G)\big),
\end{equation*}
which we call the {\em deformation space.\/} 
For applications to differential geometry, such as moduli spaces
of flat connections (gauge theory) or locally homogeneous geometric
structures, it plays a more prominent role than the representation
variety $\hpg$. Although $\Inn(G)$ preserves the algebraic structure,
$\hpgg$ will generally not admit the structure of an algebraic set.

Since the $\Inn(G)$-action on $\hpg$ absorbs the $\Inn(\pi)$-action
on $\hpg$, the outer automorphism group $\Out(\pi)$ acts on
$\hpgg$. 
By a theorem of M.\ Dehn and J.\ Nielsen (compare Nielsen~\cite{Nielsen1927}
and Stillwell~\cite{Stillwell}), $\Out(\pi)$ identifies with the
{\em mapping class group\/} 
\begin{equation*}
\Mod :=  \pi_0\big(\Diff(\Sigma)\big).
\end{equation*}
One motivation for this study is that the deformation
spaces provide natural objects upon which mapping class 
groups act~\cite{Goldman_mcg}.

\section{Abelian groups and rank one Higgs bundles}\label{sec:abelian}
The simplest groups are commutative. 
When $G$ is abelian, then the commmutators $[\alpha,\beta]=1$ and
the defining relation in \eqref{eq:presentation} is vacuous.
Thus
\begin{equation*}
\hpg \longleftrightarrow G^{2g}
\end{equation*}
Furthermore $\Inn(G)$ is trivial so
\begin{equation*}
\hpgg \longleftrightarrow G^{2g}
\end{equation*}
as well.

\subsection{Symplectic vector spaces}

Homological machinery applies. By the Hurewicz theorem and the universal
coefficient theorem,
\begin{equation*}
\hpg \cong \Hom(\pi/[\pi,\pi],G) \cong \Hom(H_1(\Sigma),G) \cong
H^1(\Sigma,G)
\end{equation*}
(or $H^1(\pi,G)$ if you prefer group cohomology). In particular
when $G= \R$, then $\hpgg$ is the real vector space
\begin{equation*}
H^1(\Sigma,\R) \cong \R^{2g}
\end{equation*}
which is naturally a {\em symplectic vector space\/} under the cup-product
pairing
\begin{equation*}
H^1(\Sigma,\R) \times H^1(\Sigma,\R) \longrightarrow
H^2(\Sigma,\R)\cong \R,
\end{equation*}
the last isomorphism corresponding to the orientation of $\Sigma$.

Similarly when $G = \C$, the representation variety and the deformation
space
\begin{equation*}
\hpgg = \hpg \longleftrightarrow H^1(\pi,\C) \cong H^1(\Sigma,\C) 
\end{equation*}
is a {\em complex-symplectic vector space,\/} that is, a complex vector
space with a complex-bilinear symplectic form.

The mapping class group action factors through the action on 
homology of $\Sigma$, or equivalently the abelianization of $\pi$, 
which is the homomorphism
\begin{equation*}
\Mod \longrightarrow \Sp(2g,\Z).
\end{equation*}

\subsection{Multiplicative characters: $G=\C^*$}

Representations $\pi \longrightarrow \C^*$ 
correspond to {\em multiplicative characters,\/} and are easily
understood using the universal covering
\begin{align*}
\C &\longrightarrow \C^* \\
z &\longmapsto \exp(2\pi i z)
\end{align*}
with kernel $\Z\subset\C$.
Such a representation corresponds to a {\em flat complex line bundle
over $\Sigma$.\/}
The deformation space $\hpg$ identifies with the quotient
\begin{equation*}
H^1(\Sigma,\C) / H^1(\Sigma,\Z).
\end{equation*}
Restricting to unit complex numbers $G = \Uo \subset \C^*$, 
identifies $\hpg$ with the $2g$-dimensional {\em torus\/}
\begin{equation*}
H^1(\Sigma,\R) / H^1(\Sigma,\Z),
\end{equation*}
the quotient of a real symplectic vector space by an integer lattice,
$\Mod$-acts on this torus by {\em symplectomorphisms.\/}

\subsection{The Jacobi variety of a Riemann surface}

The classical Abel-Jacobi theory 
(compare for example Farkas-Kra~\cite{FarkasKra}),
identifies unitary characters $\pi_1(X)\longrightarrow\Uo$ 
of the fundamental group of a Riemann surface
$X$ with topologically trivial holomorphic line bundles over $X$.
In particular $\hpg$ identifies with the {\em Jacobi variety\/}
$\Jac$.

While the basic structure of $\hpg$ is a $2g$-dimensional compact
real torus with a parallel symplectic structure, the conformal structure
on $X$ provides much stronger structure. Namely, $\Jac$ is a 
{\em principally polarized abelian variety,\/} a projective variety with the
structure of an abelian group. Indeed this extra structure, by Torelli's
theorem, is enough to recover the Riemann surface $X$.

In particular the analytic/algebraic structure on $\Jac$ is definitely
{\em not\/} invariant under the mapping class group $\Mod$. However the 
symplectic structure on $\hpg$ is independent of the conformal structure $X$
and is invariant under $\Mod$.

The complex structure on $\Jac$ is the effect of the complex structure
on the tangent bundle $TX$ (equivalent to the Hodge $\star$-operator). 
The Hodge theory of harmonic differential forms finds
unique harmonic representatives for cohomology classes, which uniquely
extend to {\em holomorphic differential forms.\/} Higgs bundle theory
is {\em nonabelian  Hodge theory\/} (Simpson~\cite{Simpson2})
in that it extends this basic technique from ordinary 1-dimensional
cohomology classes to flat connections. 

When $G = \C^*$, then $\hpg$ acquires a
complex structure $J$ coming from the complex structure on $\C^*$.
This depends only on the topology $\Sigma$, in fact just its fundamental
group $\pi$. 
Cup product provides a holomorphic symplectic
structure $\Omega$ 
on this complex manifold, giving the moduli space the structure
of a {\em complex-symplectic manifold.\/}

As for the $\Uo$-case  above, Hodge theory on the Riemann surface $X$ 
determines another complex structure
by $I$; then these two complex structures anti-commute:
\begin{equation*}
I J + J I = 0, 
\end{equation*}
generating a quaternionic action on the tangent bundle with $K:=IJ$. 
The symplectic structure arising from cup-product is not holomorphic
with respect to $I$; instead it is {\em Hermitian\/} (of Hodge type
$(1,1)$) with respect to $I$, extending the K\"ahler structure on $\Jac$.
Indeed with the structure $I$, $\hp{\C^*}$ identifies with the 
{\em cotangent bundle\/} $T^*\Jac$ with K\"ahler metric defined by
\begin{equation*}
g(X,Y) :=  \Omega(X, IY)
\end{equation*}

The triple $(\Omega,I,J)$ defines a {\em \hK structure \/} 
{\em refining\/} the complex-sym\-plectic structure.
If one thinks of a complex-symplectic structure as a $G$-structure
where $G = \Sp(2g,\C)$, then a \hK refinement is a reduction of the
structure group to the maximal compact $\Sp(2g,\C) \supset \Sp(2g)$.
The more common definition of a \hK structure 
involves the Riemannian metric $g$ which
is K\"ahlerian with respect to all three complex structures $I,J,K$;
alternatively it is characterized as a Riemannian manifold of dimension $4g$ 
with holonomy reduced to $\mathsf{Sp}(2g)\subset \mathsf{SO}(4g)$.

For a detailed exposition of the theory of rank one Higgs bundles on 
Riemann surfaces, compare Goldman-Xia~\cite{GoldmanXia}.

\section{Stable vector bundles and Higgs bundles}\label{sec:Stable}
Narasimhan and Seshadri~\cite{NarasimhanSeshadri} generalized the
Abel-Jacobi theory above to identify $\hpgg$ with a moduli space
of holomorphic objects over a Riemann surface $X$, when $G=\Un$.
(This was later extended by Ramanathan~\cite{Ramanathan} to general
compact Lie groups $G$.)

A notable new feature is that, unlike line bundles, not every
topologically trivial holomorphic rank $n$ vector bundle arises from
from a representation $\pi\longrightarrow\Un$. Furthermore 
equivalence classes of all holomorphic $\C^n$-bundles does not form an
algebraic set.

Narasimhan and Seshadri define a degree zero holomorphic $\C^n$-bundle
$V$ over $X$ to be {\em stable\/} (respectively {\em semistable\/})
if and only if 
every holomorphic vector subbundle of $V$ has
negative (respectively nonpositive) degree.  Then a holomorphic vector
bundle arising from a unitary representation $\rho$ is semistable, and
the bundle is stable 
if and only if 
the representation is
irreducible.  Furthermore every such semistable bundle arises from a
unitary representation. Narasimhan and Seshadri show the moduli space
$\NS$ of semistable bundles of degree 0 and rank $n$
over $X$ is naturally a projective variety, thus
defining such a structure on $\hpgg$. The K\"ahler structure
depends heavily on the Riemann surface $X$, although the symplectic
structure depends only on the topology $\Sigma$.

%%%%
It is useful to extend the notions of stability to bundles which
may not have degree $0$. In particular we would like stability to
be preserved by tensor product with holomorphic line bundles. 
Define a holomorphic vector bundle $V$ to be {\em stable\/} 
if every holomorphic subbundle $W\subset V$ satisfies the
inequality
\begin{equation*}
\frac{\deg(W)}{\rank(W)} <  \frac{\deg(V)}{\rank(V)}.
\end{equation*}
Semistability is defined similarly by replacing the strict
inequality by a weak inequality.
%%%%

In trying to extend this correspondence to the complexification $G =
\GLnC$ of $\Un$, one might consider the {\em cotangent bundle\/}
$T^*\NS$ of the Narasimhan-Seshadri moduli space, and relate it to
representations $\pi\to\GLnC$. In particular since
cotangent bundles of K\"ahler manifolds tend to be 
hyper-K\"ahler, relating
$\hpgg$ to $T^\NS$ might lead to a \hK~geometry on $\hpgg$.

Thus a neighborhood of the $\Un$-representations in the space of
$\GLnC$ corresponds to a neighborhood of the zero-section of
$T^*\NS$. In turn, elements in this neighborhood identify with pairs
$(V,\Phi)$ where $V$ is a semistable holomorphic vector bundle and
$\Phi$ is a tangent covector to $V$ in the space of holomorphic vector
bundles. Such a tangent covector is with a 
{\em Higgs field,\/} by definition, an $\End(V)$-valued holomorphic $1$-form 
on $X$.

Although one can define a \hK~structure on the moduli space of such pairs,
the \hK~metric is incomplete and not all irreducible linear representations
arise. To rectify this problem, one must consider Higgs fields on possibly
unstable vector bundles.

Following Hitchin~\cite{Hitchin_selfduality} and
Simpson~\cite{Simpson1}, define a {\em Higgs pair\/} to be a pair
$(V,\Phi)$ where $V$ is a (not ncessarily semistable) holomorphic vector 
bundle and the Higgs field $\Phi$ a $\End(V)$-valued holomorphic $1$-form. 
Define $(V,\Phi)$ to be 
{\em stable\/} if and only if 
for all  $\Phi$-invariant holomorphic subbundles $W\subset V$,
\begin{equation*}
\frac{\deg(W)}{\rank(W)} <  \frac{\deg(V)}{\rank(V)}.
\end{equation*}
The Higgs bundle $(V,\Phi)$ is {\em polystable\/} 
if and only if
$(V,\Phi) = \bigoplus_{i=1}^l (V_i,\Phi_i)$
where each summand $(V_i,\Phi_i)$ is stable and 
\begin{equation*}
\frac{\deg(V_i)}{\rank(V_i)} =  \frac{\deg(V)}{\rank(V)}
\end{equation*}
for $i=1,\dots,l$.

The following basic result follows from 
Hitchin~\cite{Hitchin_selfduality}, Simpson~\cite{Simpson1},
with a key ingredient (the {\em harmonic metric\/}) 
supplied by Corlette~\cite{Corlette_88} and Donaldson~\cite{Donaldson}:

\begin{thm*}
The following natural bijections exist between equivalences classes:

\begin{equation*}
\left\{ \begin{aligned}
& \text{\rm Stable Higgs pairs}\\ & 
(V,\Phi) \text{\rm ~over~} \Sigma
\end{aligned}
\right\}  
\longleftrightarrow 
\left\{ \begin{aligned}
& \text{\rm Irreducible representations }\\ & 
\pi_1(\Sigma)\xrightarrow{\rho} \GLnC
\end{aligned}
\right\}
\end{equation*}

\begin{equation*}
\left\{ \begin{aligned}
& \text{\rm Polystable Higgs pairs}\\ & 
(V,\Phi) \text{\rm ~over~} \Sigma
\end{aligned}
\right\}  
\longleftrightarrow 
\left\{ \begin{aligned}
& \text{\rm Reductive representations }\\ & 
\pi_1(\Sigma)\xrightarrow{\rho} \GLnC
\end{aligned}
\right\}
\end{equation*}
\end{thm*}
\noindent
When the Higgs field $\Phi=0$, 
this is just the Narasimhan-Seshadri theorem, 
identifying stable holomorphic vector bundles 
with irreducible $\U(n)$-representations.
Allowing the Higgs field $\Phi$ to be nonzero, even when $V$ is unstable,
leads to a rich new class of examples, which can now be treated using
the techniques of Geometric Invariant Theory.

\section{Hyperbolic geometry: $G = \PSLtR$}

Another important class of surface group representations are 
{\em Fuchsian representations,\/} 
which arise by endowing  $\Sigma$ with the
local geometry of {\em hyperbolic space $\Ht$.\/} Here $G$ is the
group of orientation-preserving isometries $\Isom(\Ht)$, which, using
Poincar\'e's {\em upper half-space model,\/} identifies with $\PSLtR$.
Fuchsian representations are characterized in many different
equivalent ways; in particular a representation $\pi\xrightarrow{\rho}
G = \PSLtR$ is Fuchsian 
if and only if
it is a {\em discrete
embedding,\/} that is, $\rho$ embeds $\pi$ isomorphicly onto a
discrete subgroup of $G$.

\subsection{Geometric structures}

Let $\Ht$ be the hyperbolic plane with a fixed orientation 
and $G \cong \Isom(\Ht) \cong\PSLtR$ 
its group of orientation-preserving isometries.
A {\em hyperbolic structure\/} on a topological surface $\Sigma$ is defined
by a coordinate atlas $\{ (U_\alpha,\psi_\alpha) \}_{\alpha\in A}$ where
\begin{itemize}
\item
The collection $\{U_\alpha\}_{\alpha\in A}$ of coordinate patches 
{\em covers\/} $\Sigma$ (for some index set $A$);
\item
Each coordinate {\em chart\/} $\psi_\alpha$ is an orientation-preserving 
homeomorphism of the coordinate patch $U_\alpha$ onto an open subset
$\psi_\alpha(U_\alpha)\subset\Ht$.
\item
For each connected component $C\subset U_\alpha\cap U_\beta$, there is
(necessarily unique) $g_{C,\alpha,\beta} \in G$ such that
\begin{equation*}
\psi_\alpha |_C  = 
g_{C,\alpha,\beta} \circ \psi_\beta |_C.
\end{equation*}
\end{itemize}
\noindent

The resulting {\em local\/} hyperbolic geometry defined on the patches
by the coordinate charts is independent of the charts, and extends to
a global structure on $\Sigma$. The surface $\Sigma$ with this refined
structure of local hyperbolic geometry, will be called a {\em
hyperbolic surface\/} and denoted by $M$.
Such a structure is equivalent to a Riemannian metric
of constant curvature $-1$. The equivalence follows from two basic facts: 
\begin{itemize}
\item 
Any two Riemannian manifolds 
of curvature $-1$ are locally isometric;,
\item
A local isometry from a connected subdomain of $\Ht$ 
extends globally to an isometry of $\Ht$.
\end{itemize}
\noindent
Suppose $M_1,M_2$ are two hyperbolic surfaces. 
Define a {\em morphism\/} $M_1\xrightarrow{\phi} M_2$ 
as a map $\phi$, which, in the preferred local
coordinates of $M_1$ and $M_2$, is defined by isometries in
$G$. Necessarily a morphism is a local isometry of Riemannian
manifolds. Furthermore, if $M$ is a hyperbolic surface and
$\Sigma\xrightarrow{f} M$ is a local homeomorphism, there exists a
hyperbolic structure on $\Sigma$ for which $f$ is a morphism.
In particular every covering space of a hyperbolic
surface is a hyperbolic surface.

In more traditional terms, a morphism of hyperbolic surfaces is just a
local isometry.

\subsection{Relation to the fundamental group}

While the definitions involving coordinate atlases or Riemannian metrics
have certain advantages, another point of view underscores
the role of the fundamental group.

Let $M$ be a hyperbolic surface.
Choose a universal covering space $\tM\to M$ and give $\tM$
the unique hyperbolic structure for which %the covering projection 
$\tM\to M$ is a local isometry. 
Then there exists a {\em developing map\/} $\tM \xrightarrow{\dev_M} \Ht$,
a local isometry,  which induces the hyperbolic structure on $\tM$ from that
of $\Ht$.
The group $\pi_1(M)$ of deck transformations of
$\tM\to M$ acts on $\Ht$ by isometries and $\dev$ is equivariant
respecting this action: for all $\gamma\in\pi_1(M)$, the diagram
\begin{equation*}
\begin{CD}
\tM @>{\dev_M}>>  \Ht \\
@V{\gamma}VV    @VV{\rho(\gamma)}V \\
\tM @>>{\dev_M}>   \Ht.
\end{CD}
\end{equation*}
commutes. The correspondence $\gamma \longmapsto \rho(\gamma)$ is 
a homomorphism, 
\begin{equation*}
\pi_1(M)\xrightarrow{\hol_M} \Isom(\Ht),  
\end{equation*}
the  {\em holonomy representation\/} of the hyperbolic surface $M$.
The pair $(\dev_M,\hol_M)$ is unique up to the $G$-action defined by
\begin{equation*}
(\dev_M,\hol_M) \stackrel{g}\longmapsto 
(g\circ\dev_M,\Inn(g)\circ\hol_M)
\end{equation*}
for $g\in\Isom(\Ht)$.

If the hyperbolic structure is {\em complete,\/} that is, the
Riemannian metric is geodesically complete, then the developing map
is a {\em global isometry\/} $\tM \approx \Ht$. In that case
the $\pi$-action on $\Ht$ defined by the  holonomy representation $\rho$ 
is equivalent to the action by deck transformations.
Thus $\rho$
defines a proper free action of $\pi$ on $\Ht$ by 
isometries. Conversely if $\rho$ defines a proper free 
isometric $\pi$-action, then the quotient 
\begin{equation*}
M :=\Ht/\rho(\pi) 
\end{equation*}
is a complete hyperbolic manifold with a preferred isomorphism 
\begin{equation*}
\pi_1(\Sigma) \xrightarrow{\rho} \rho(\pi) \subset G.
\end{equation*}
This isomorphism (called a {\em marking\/}) 
determines a preferred homotopy class of homotopy equivalences
\begin{equation*}
\Sigma \longrightarrow M. 
\end{equation*}

\subsection{Examples of hyperbolic structures}
We now give three examples of surface group representations in $\PSLtR$.  
The first example is Fuchsian and corresponds to a hyperbolic structure on
a surface of genus two. The second example is not Fuchsian, but
corresponds to a hyperbolic structure with a single {\em branch
point,\/} that is a point with local coordinate given by a branched
conformal mapping $z \longmapsto z^k$ where $k \ge 1$.  (The
nonsingular case corresponds to $k=1$.) In our example $k=2$ and the
singular point has a neighborhood isometric to a hyperbolic cone of
cone angle $4\pi$. 

\subsubsection{A Fuchsian example}\label{sec:45}
Here is a simple example of a hyperbolic surface of genus two.
Figure~\ref{fig:octagon} depicts a topological construction for
a genus two surface $\Sigma$. Realizing this topological construction
in hyperbolic geometry gives $\Sigma$ a local hyperbolic geometry as follows.
Take a regular octagon $P$ with angles $\pi/4$. Label the sides as 
\begin{equation*}
A_1^-,B_1^-,A_1^+,B_1^+, 
A_2^-,B_2^-,A_2^+,B_2^+ 
\end{equation*}
$a_i$ pairs $B_i^-$ to $B_i^+$
and $b_i$ pairs $A_i^-$ to $A_i^+$ respectively.

Pair the sides by 
\begin{equation*}
a_1,b_1,a_2,b_2\in\PSLtR 
\end{equation*}
according to the pattern described in Figure~\ref{fig:octagon}.
Given any two oriented geodesic segments in $\Ht$ of equal length,
a unique orientation-preserving isometry maps one to the other. 
Since the polygon is regular, one can realize all four identifications
in $\Isom(\Ht)$.

The quotient (compare Figure~\ref{fig:genustwo}) contains three types
of points: 
\begin{itemize}
\item 
A point in the open 2-cell has a coordinate
chart which is the embedding $P\hookrightarrow\Ht$.
\item
A point on the interior of an edge has a half-disc neighborhood, which
together with the half-disc neighborhood of its part, gives a coordinate
chart for the corresponding point in the quotient. 
\item
Around the single vertex in the quotient is a cone of angle 
\begin{equation*}
8 (\pi/4) = 2\pi,
\end{equation*}
a disc in the hyperbolic plane. 
\end{itemize}
The resulting identification space is a hyperbolic surface of genus
$g=2$.
The above isometries satisfying the defining relation for $\pi_1(\Sigma)$:
\begin{equation*}
a_1 b_1 a_1^{-1} b_1^{-1}   a_2 b_2 a_2^{-1} b_2^{-1}  = 1
\end{equation*}
and define a Fuchsian representation 
\begin{equation*}
\pi_1(\Sigma) \xrightarrow{\rho} \PSLtR.
\end{equation*}

\begin{figure}[th]
\centerline{\epsfig{figure=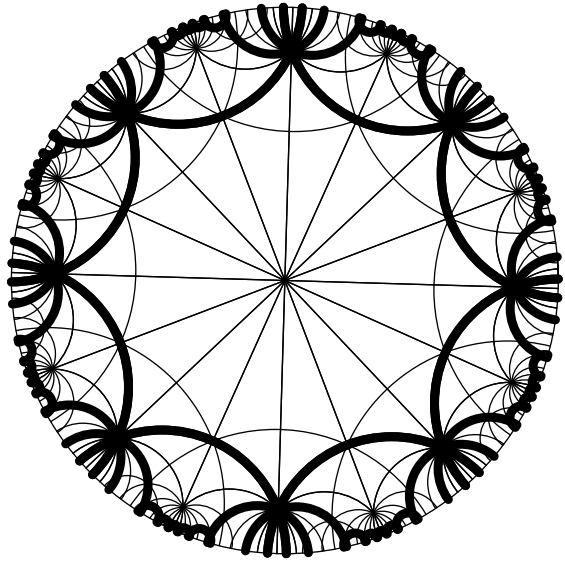,height=8cm}}
\caption{ 
A regular octagon with vertex angles $\pi/4$ can be realized in the
tiling of $\Ht$ by triangles with angles $\pi/2,\pi/4,\pi/8$. 
The identifications depicted in Figure~\ref{fig:octagon} are realized
by orientation-preserving isometries. The eight angles of $\pi/4$ fit together 
to form a cone of angle $2\pi$, forming a coordinate chart for a hyperbolic
structure around that point.
}
\label{fig:45tiling}
\end{figure}

\subsubsection{A branched hyperbolic structure}\label{sec:90}
We can modify the preceding example to include a singular structure,
again on a surface of genus two.
Take a regular {\em right-angled \/} octagon. Again,
labeling the sides as before, side pairings $a_1,b_1,a_2,b_2$
exist. Now $8$ right angles compose a neighborhood of the vertex in
the quotient space. The quotient space is a hyperbolic structure
with one singularity of cone angle
$4 \pi = 8 (\pi/2)$.
Since the product of the identification mappings
\begin{equation*}
a_1 b_1 a_1^{-1} b_1^{-1}   a_2 b_2 a_2^{-1} b_2^{-1} 
\end{equation*}
is rotation through $4\pi$ (the identity),
the holonomy representation $\hat\rho$ 
of the nonsingular hyperbolic surface 
$\Sigma\setminus \{p\}$ extends:
\begin{equation*}
\xymatrix{
& \pi_1\big(S \setminus \{p\}\big) 
\ar@{>>}[d] \ar[dr]^{\hat\rho} \\
  & \pi_1(\Sigma) \ar@{-->}[r]_{\rho} & \PSLtR }
\end{equation*}
% and $\Euler(\rho) = -1$. 
Although $\rho(\pi)$ is discrete,
% (a Fuchsian group representing a torus with a branch 
% point of cone angle $\pi$),
$\rho$ is not injective.

\begin{figure}[ht]
\centerline{\epsfig{figure=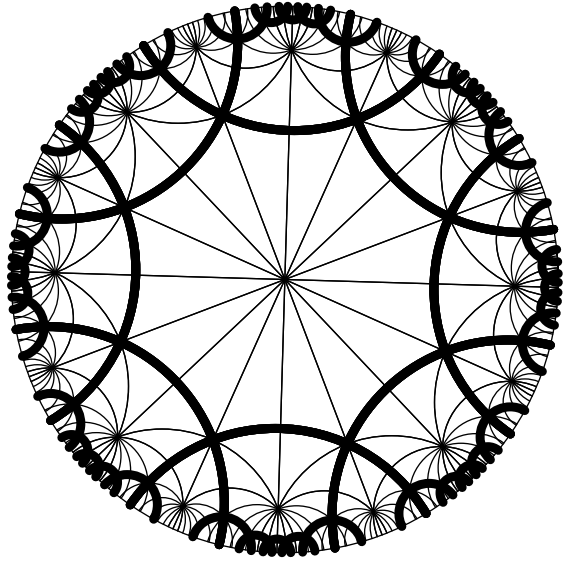,height=8cm}}
\caption{ 
A regular right-angled octagon can also be realized in the
tiling of $\Ht$ by triangles with angles $\pi/2,\pi/4,\pi/8$. 
The identifications depicted in Figure~\ref{fig:octagon} are realized
by orientation-preserving isometries. The eight angles of $\pi/2$ fit together 
to form a cone of angle $4\pi$, forming a coordinate chart for a 
singular hyperbolic structure, branched at one point.
}
\label{fig:90tiling}
\end{figure}

\subsubsection
{A representation with no branched structures}\label{sec:nobranch}

\begin{figure}[ht]
\centerline{\epsfig{figure=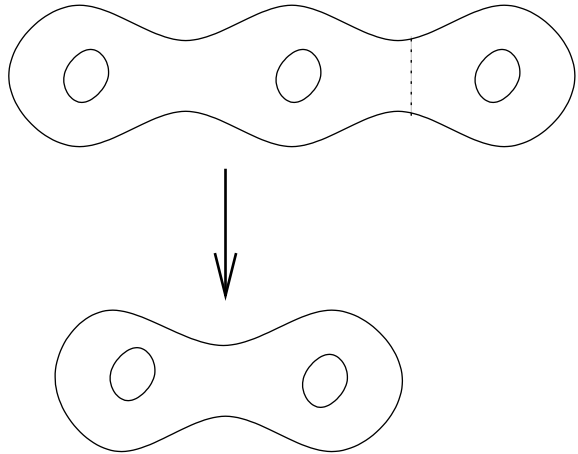,height=6cm}}
\caption{A degree one map from a genus 3 surface to a genus 2
surface which collapses a handle. Such a map is not homotopic
to a smooth map with branch point singularities (such as a holomorphic
map).}
\label{fig:genusthree}
\end{figure}

Consider a degree-one map $f$ from a genus three surface $\Sigma$ to
a genus two hyperbolic surface $M$, depicted in Figure~\ref{fig:genusthree}.
Let $\pi_1(M)\xrightarrow{\mu} G$ 
denote the holonomy representation of $M$ and consider the composition
\begin{equation*}
\pi = \pi_1(\Sigma) \xrightarrow{f_*} \pi_1(M) \stackrel{\mu}\hookrightarrow G.
\end{equation*}
Then a branched hyperbolic structure with holonomy $\mu\circ f_*$
corresponds to a mapping with branch singularities
\begin{equation*}
\Sigma \xrightarrow{F} \Ht/\mathsf{Image}(\mu\circ f_*) = M.
\end{equation*}
inducing the homomorphism
\begin{equation*}
\pi = \pi_1(\Sigma) \xrightarrow{f_*} \pi_1(M).
\end{equation*}
In particular $F\simeq f$. However, since $\deg(f) = 1$, any mapping with
only branch point singularities of degree one must be a homeomorphism,
a contradiction.

\section{Moduli of hyperbolic structures and representations}
To understand ``different'' geometric structures on the ``same''
surface, one introduces {\em markings.\/} 
Fix a topological type
$\Sigma$ and let the geometry $M$ vary.  The 
fundamental group $\pi = \pi_1(\Sigma)$ is also fixed, 
and each marked structure determines a well-defined equivalence class
in $\hpgg$.
Changing the marking corressponds to the action of the mapping
class group $\Mod = \Out(\pi)$ on $\hpgg$. 
{\em Unmarked structures\/} 
correspond to the orbits of the $\Mod$-action.

\subsection{Deformation spaces of geometric structures}

A {\em marked hyperbolic structure on $\Sigma$\/} is defined as
a pair $(M,f)$ where $M$ is a hyperbolic surface and $f$
is a homotopy equivalence $\Sigma \to M$.
Two marked hyperbolic structures 
\begin{equation*}
\Sigma \xrightarrow{f} M, \qquad
\Sigma \xrightarrow{f'} M'.  
\end{equation*}
are {\em equivalent\/} if and only if there exists an isometry
$M\xrightarrow{\phi}M'$ such that 
 \begin{equation*}
 \xymatrix{
& M \ar@{-->}[d]^{\phi}\\
\Sigma \ar@{>}[ur]^{f}
\ar@{>}[r]_{f'} & M' 
}
\end{equation*}
homotopy-commutes, that is, $\phi\circ f\simeq f'$.  
The {\em Fricke space\/} $\Fricke(\Sigma)$
of $\Sigma$ is the space of all such equivalence classes of
marked hyperbolic structures on $\Sigma$.  
(Bers-Gardiner~\cite{BersGardiner}.)
The Fricke space is diffeomorphic to $\R^{6g-6}$.
The theory of moduli of hyperbolic structures on
surfaces goes back at least to 
Fricke and Klein~\cite{FrickeKlein}.

The Teichm\"uller space $\Teich$ of $\Sigma$ is defined similarly, as
the space of equivalence classes of {\em marked conformal
structures\/} on $\Sigma$, that is, pairs $(X,f)$ where $X$ is a
Riemann surface and $\Sigma\xrightarrow{f} X$ is a homotopy
equivalence.  Teichm\"uller used quasiconformal mappings
to parametrize $\Teich$ by elements of a vector space, define a metric
on $\Teich$ and prove analytically that $\Teich$ is a cell.  
Using these ideas, Ahlfors~\cite{Ahlfors})
proved $\Teich$ is naturally a complex manifold.

Sine a hyperbolic structure is a Riemannian
metric, every hyperbolic structure 
has an underlying conformal structure.
The {\em uniformization theorem\/} asserts that if
$\chi(\Sigma) < 0$, then every conformal structure
on $\Sigma$ underlies a unique hyperbolic structure.
The resulting identification of conformal and hyperbolic
structures identifies $\Teich$ with $\Fricke(\Sigma)$.
As discussed below, $\Fricke(\Sigma)$ identifies
with an open subset of $\Hom(\pi,\PSLtR)/\PSLtR$ which
has no apparent complex structure. Thus the complex
structure on $\Teich$ is more mysterious when
$\Teich$ is viewed as a space of hyperbolic structures.
For a readable survey of classical Teichm\"uller theory see
Bers~\cite{Bers_survey}.

\subsection{Fuchsian components of $\hpgg$}

To every equivalence class of marked hyperbolic structures is associated
a well-defined element
\begin{equation*}
[\rho] \in \hpgg.
\end{equation*}
A representation $\pi\xrightarrow{\rho} G$ is {\em Fuchsian \/}
if and only if
it arises as the holonomy of a hyperbolic structure on $\Sigma$.
Equivalently, it satisfies the three conditions:
\begin{itemize}
\item $\rho$ is injective;
\item Its image $\rho(\pi)$ is a discrete subgroup of $G$;
\item The quotient $G/\rho(\pi)$  is compact.
\end{itemize}
The first condition asserts that $\rho$ is an {\em embedding,\/}
and the second two conditions assert that $\rho(\pi)$ is a
{\em cocompact lattice.\/} Under our assumption $\partial\Sigma=\emptyset$, 
the third condition (compactness of $G/\rho(\pi)$) follows from the
first two. In general, we say that $\rho$ is a {\em discrete embedding\/}
(or {\em discrete and faithful\/}) 
if $\rho$ is an embedding with discrete image (the first two conditions).

\begin{thm}
Let $G=\mathsf{Isom}(\Ht) = \PGL(2,\R)$ and $\Sigma$ a closed connected
surface with $\chi(\Sigma)<0$.  Fricke space, the subset of $\hpgg$
consisting of $G$-equivalence classes of Fuchsian representations, is
a connected component of $\hpgg$.
\end{thm}

This result follows from three facts: 
\begin{itemize}
\item Openness of Fricke space
(Weil~\cite{Weil}), 
\item Closedness of Fricke space
(Chuckrow~\cite{Chuckrow}), 
\item Connectedness of Fricke space 
\end{itemize}

Chuckrow's theorem is a special case of a consequence of 
the Kazhdan-Margulis uniform discreteness 
(compare Raghunathan~\cite{Raghunathan} 
and Goldman-Millson~\cite{Goldman_Millson1}).
These ideas go back to Bieberbach and Zassenhaus
in connection with the classification of Euclidan crystallographic
groups.
Uniform discretess applies under very
general hypotheses, to show that discrete embeddings form a closed
subset of the representation variety.
For the proof of connectedness, see,
for example, Jost~\cite{Jost},\S 4.3, Buser~\cite{Buser}, 
Thurston~\cite{Thurston}
or Ratcliffe~\cite{Ratcliffe} 
for elementary proofs using Fenchel-Nielsen
coordinates).
Connectedness also follows from the uniformization
theorem, together with the identification of Teichm\"uller space 
$\Teich$ as a cell.

When $G=\Isom(\Ht)=\PSL(2,\R)$, the situation slightly 
complicates, due to the choice of orientation. 
Assume $\Sigma$ is orientable,
and orient it. Orient $\Ht$ as well. 
Let $\Sigma\xrightarrow{f}M$ be a marked hyperbolic structure on $\Sigma$.
The orientation of $M$ induces an orientation of $\tM$ which is invariant
under $\pi_1(M)$. However, the developing map $\dev_M$ may or not preserve
the (arbitrarily chosen) orientations of $\tM$ and $\Ht$. Accordingly
$\Isom(\Ht)$-equivalence classes of Fuchsian representations in $G$ 
fall into two classes, which we call {\em orientation-preserving\/} and
{\em orientation-reversing\/} respectively. These two classes are
interchanged by inner automorphisms of orientation-reversing isometries
of $\Ht$.

\begin{thm}
Let $G=\Isom(\Ht) = \PSLtR$ and $\Sigma$ a closed connected oriented
surface with $\chi(\Sigma)<0$.  The set of $G$-equivalence classes of
Fuchsian representations forms two connected connected components of
$\hpgg$. One component corresponds to orientation-preserving Fuchsian 
representations and the other to orientation-reversing Fuchsian 
representations.
\end{thm}

\subsection{Characteristic classes and maximal representations}
Characteristic class\-es of flat bundles determine invariants of
representations. In the simplest cases (when $G$ is compact or
reductive complex), these determine the connected components of $\hpg$.

\subsubsection{The Euler class and components of $\Hom(\pi,\PSLtR)$}
The components of $\hpg$ were determined in \cite{Goldman_components}
using an invariant derived from the Euler class of the oriented $\Ht$-bundle
\begin{equation*}
\xymatrix{
\Ht
 \ar@{>}[r] 
& (\Ht)_\rho  \ar@{>}[d] \\ & \Sigma
}
\end{equation*}
associated to a representation $\pi\xrightarrow{\rho}\PSLtR$ as follows.
The total space is the quotient
\begin{equation*}
(\Ht)_\rho  := (\tilde\Sigma\times \Ht)/\pi
\end{equation*}
where $\pi$ acts diagonally on $\tilde\Sigma\times \Ht$ by deck transformations
on $\tilde\Sigma$ and via $\rho$ on $\Ht$.
Isomorphism classes of oriented $\Ht$-bundles over $\Sigma$ are determined 
by the {\em Euler class,\/} which lives in $H^2(\Sigma,\Z)$. The orientation
of $\Sigma$ identifies this cohomology group with $\Z$.
The resulting map
\begin{equation*}
\hpg \xrightarrow{\Euler} H^2(S;\Z) \cong \Z
\end{equation*}
satisfies
\begin{equation}\label{eq:MilnorWood}
\vert \Euler(\rho) \vert \le \vert\chi(S)\vert = 2 g -2 .
\end{equation}
(Milnor~\cite{Milnor} and  Wood~\cite{Wood}).
Call a representation {\em maximal\/} if equality holds in
in \eqref{eq:MilnorWood}, that is, 
$\Euler(\rho) = \pm\chi(\Sigma)$: 

The following converse was proved in Goldman~\cite{Goldman_thesis} 
(compare also \cite{Goldman_components} and \cite{Hitchin_selfduality}).
\begin{thm}
$\rho$ is maximal if and only if $\rho$ is Fuchsian.
\end{thm}
Suppose $M$ is a branched hyperbolic surface with branch points 
$p_1,\dots, p_l$
where $p_i$ is branched of order $k_i$, where each $k_i$ is a positive
integer. In other words, each $p_i$
has a neighborhood which is a hyperbolic cone of cone angle $2\pi k_i$.
Consider a marking $\Sigma \to M$, determining a holonomy 
representation $\rho$. Then
\begin{equation*}
\Euler(\rho) = \chi(\Sigma) + \sum_{i=1}^l k_i.
\end{equation*}
Consider the two examples for genus two surfaces.
\begin{itemize}
\item 
The first (Fuchsian) example
(\S\ref{sec:45})
arising from a regular octagon with $\pi/4$ angles, has Euler class
$-2 = \chi(\sigma)$. 
\item In the second example (\S\ref{sec:90}), 
the structure is branched at one point, so that  $l=k_1=1$ and
the Euler class equals $-1 = \chi(\Sigma) + 1$. 
\end{itemize}
\subsection{Quasi-Fuchsian representations: $G=\PSLtC$}
When the representation 
\begin{equation*}
\pi \longrightarrow \PSLtR \hookrightarrow \PSLtC 
\end{equation*}
is deformed inside $\PSLtC$,
the action on $\cpo$ is topologically conjugate
to the original Fuchsian action. 
Furthermore there exists a H\"older $\rho$-equivariant embedding 
$S^1 \hookrightarrow \cpo,  $
whose image $\Lambda$ has Hausdorff
dimension $>1$, --- unless the deformation is still Fuchsian.
The space of such representations is the {\em quasi-Fuchsian space\/} $\QF$.
By Bers~\cite{Bers}, $\QF$ naturally identifies with 
\begin{equation*}
\Teich\times\overline{\Teich} \approx \R^{12g-12}.
\end{equation*}
%%%%
Bers's correspondence is the following. 
The action of $\rho$ on the complement $\cpo\setminus\Lambda$
is properly discontinuous, and the quotient 
\begin{equation*}
\big(\cpo\setminus\Lambda \big)/\rho(\pi)
\end{equation*}
consists of two Riemann surfaces, 
each with a canonical marking determined by $\rho$.
Furthermore these surfaces possess opposite orientations,
so the pair of marked conformal structures determine a point
in $\Teich\times\overline{\Teich}$.
%%%%
Bonahon~\cite{Bonahon} and Thurston proved that the 
closure of $\QF$ in $\hpgg$ equals the 
set of equivalence classes of discrete embeddings 
The frontier $\partial\QF\subset\hpgg$ 
is nonrectifiable, and is near non-discrete representations.

However, the two connected components of $\hpgg$ are 
distinguished by the characteristic 
class (related to the second Stiefel-Whitney class $w_2$) which 
detects whether a representation in $\PSLtC$ lifts to the double covering
$\SLtC \to \PSLtC$ (Goldman~\cite{Goldman_components}).
Contrast this situation with $\PSLtR$ where the discrete
embeddings form connected components, characterized by maximality.

\begin{figure}[ht]
\centerline{\epsfig{figure=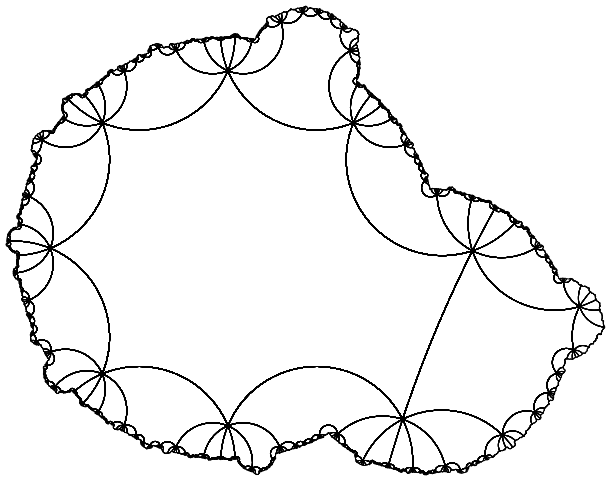,height=8cm}}
\caption{
A quasi-Fuchsian subgroup of $\PSLtC$ obtained by deforming the
genus two surface with a fundamental domain the regular octagon
with $\pi/4$ angles in $\cpo$. The limit set is a nonrectifiable
Jordan curve, but the new action of $\pi_1(\Sigma)$ is topologically
conjugate to the original Fuchsian action.
}
\label{fig:quasifuchsian}
\end{figure}

\subsubsection{Higher rank Hermitian spaces: the Toledo invariant}
\label{sec:Toledo}

Domingo Toledo~\cite{Toledo_89} generalized
the Euler class of flat $\PSLtR$-bundles
to flat $G$-bundles, where $G$ is the autormorphism group of
a Hermitian symmetric space $X$ of noncompact type.

Let $\pi \xrightarrow{\rho} G$ be a representation
and let 
\begin{equation*}
\xymatrix{
X
 \ar@{>}[r] 
& (X)_\rho  \ar@{>}[d] \\ & \Sigma
}
\end{equation*}
be the corresponding flat $(G,X)$-bundle over $\Sigma$.
Then the $G$-invariant K\"ahler form $\omega$ on $X$
defines a closed exterior $2$-form $\omega_\rho$ on the total
space $(X)_\rho$. 
Let $\Sigma\xrightarrow{s}(X)_\rho$ be a smooth section.
Then the integral
\begin{equation*}
\int_\Sigma s^*\omega_\rho
\end{equation*}
is independent of $s$, depends continuously on $\rho$ and,
after suitable normalization, assumes integer values.
The resulting {\em Toledo invariant\/}
\begin{equation*}
\hpg \xrightarrow{\tau} \Z 
\end{equation*}
satisfies
\begin{equation*}
\vert \tau(\rho)\vert \le (2g-2) \rank_\R(G).
\end{equation*}
(Domic-Toledo~\cite{DomicToledo},
Clerc-\O rsted~\cite{Clerc_Orsted}).
Define $\rho$ to be {\em maximal\/} 
if and only if
\begin{equation*}
\vert \tau(\rho)\vert = (2g-2) \rank_\R(G). 
\end{equation*}
\begin{thm}[Toledo~\cite{Toledo_89}]
$\pi\xrightarrow{\rho}\Uno $
is maximal 
if and only if
$\rho$ is a discrete embedding preserving a complex geodesic, that is,
$\rho$ is conjugate to a representation with
\begin{equation*}
\rho(\pi) \subset \U(1,1) \times \U(n-1).
\end{equation*}
\end{thm}
This rigidity has a curious consequence for the local geometry of the
deformation space. 
Let $G := \Uno$ and 
\begin{equation*}
G_0 = \Uoo\times\mathsf{U}(n-1) \subset G. 
\end{equation*}
Then, in an appropriate sense,
\begin{equation*}
\dim  \hpgg = 2g + (2g-2)\big((n+1)^2 -1\big) = (2g-2) (n+1)^2 + 2
\end{equation*}
but Toledo's rigidity result implies that the component
of maximal representations has strictly lower dimension:
\begin{equation*}
\dim  \Hom(\pi,G_0)/G_0 = 4g + (2g-2) 3 + (2g-2)((n-1)^2 - 1)
 \end{equation*}
with codimension
\begin{equation*}
8(n-1)(g-1) - 2.
\end{equation*}

\subsection{Teichm\"uller space: marked conformal structures}

The {\em Teichm\"uller space\/} $\Teich$ of $\Sigma$ is the deformation space
of marked conformal structures on $\Sigma$. 
% In this subject it is 
% convenient, although inaccurate, to call the deformation space
% of marked hyperbolic structures {\em Teichm\"uller space.\/}

A {\em marked conformal structure\/} on $\Sigma$ is 
a pair $(X,f)$ where $X$ is a Riemann surface and $f$
is a homotopy equivalence $\Sigma \to X$.
Marked conformal structures 
\begin{equation*}
\Sigma \xrightarrow{f} X, \qquad
\Sigma \xrightarrow{f'} X'.  
\end{equation*}
are {\em equivalent\/} if and only if there exists a biholomorphism
$X\xrightarrow{\phi}X'$ such that 
 \begin{equation*}
 \xymatrix{
& X \ar@{-->}[d]^{\phi}\\
 \Sigma \ar@{>}[ur]^{f}
\ar@{>}[r]_{f'} & X' 
}
\end{equation*}
homotopy-commutes.

\begin{thm}[Uniformization]
Let $X$ be a Riemann surface with $\chi(X)<0$. Then there exists a unique
hyperbolic metric whose underlying conformal structure agrees with $X$.
\end{thm}
Since every hyperbolic structure possesses an underlying conformal structure,
Fricke space $\Fricke(\Sigma)$  maps to Teichm\"uller space $\Teich$.
By the uniformization theorem, 
$\Fricke(\Sigma)\to\Teich$.
%this map 
is an isomorphism.
It is both common and tempting to confuse these two deformation spaces.
In the present context, however, it seems best to distinguish between
the representation/hyperbolic structure and the conformal structure.

For example, each Fuchsian representation determines a marked
hyperbolic structure, and hence an underlying marked conformal
structure.  An equivalence class of Fuchsian representations thus
determines a special point in Teichm\"uller space. This contrasts
sharply with other representations which do {\em not\/} generally pick
out a preferred point in $\Teich$.  
This preferred point can be characterized as the unique minimum
of an energy function on Teichm\"uller space.

% In fact, part of the power of the 
% theory is that it applies to an arbitrary conformal structure, that
% is, an {\em arbitrary Riemann surface $X$ diffeomorphic to
% $\Sigma$.\/} The main result is an identification of the Fricke space
% $\Fricke(\Sigma)$ with the complex vector space $H^0(X,K_X^2)$.

The construction, due to Tromba~\cite{Tromba}, is as follows. 
Given a hyperbolic surface $M$ and a homotopy equivalence 
$X \xrightarrow{f} M$, then by Eels-Sampson~\cite{EelsSampson}
a unique {\em harmonic map\/} $X\xrightarrow{F} M$ exists homotopic
to $f$. The harmonic map is conformal 
if and only if 
$M$ is the
uniformization of $X$. In general the nonconformality is detected by
the {\em Hopf differential\/} $\Hopf(F)\in H^0(X,K_X^2)$, defined as the 
$(2,0)$ part of the pullback by $F$ of the complexified Riemannian metric
on $M$. The resulting mapping
\begin{align*}
\Fricke(X) &\longrightarrow H^0(X,K_X^2) \\
(f,M) &\longmapsto  \Hopf(F)
\end{align*}
is a diffeomorphism.

Fixing $M$ and letting the marked complex structure $(f,X)$ vary over 
$\Teich$ yields an interesting invariant discussed in Tromba~\cite{Tromba},
and extended in Goldman-Wentworth~\cite{Goldman_Wentworth} 
and Labourie~\cite{Labourie_energy}. The {\em energy of the harmonic map\/}
$F = F(f,X,M)$ is a real-valued function on $\Teich$. In the present context
it is the square of the $L^2$-norm of $\Hopf(F)$. 

\begin{thm}[Tromba]
The resulting function $\Teich\to \R$ 
% (which depends on the marked hyperbolic surface $\Sigma\xrightarrow{f}M$,
% or equivalently its Fuchsian holonomy representation $\rho$) 
is proper,
convex, and possesses a unique minimum at the uniformaization structure $X$.
\end{thm}
For more applications of this energy function to surface group
representations, compare Goldman-Wentworth~\cite{Goldman_Wentworth}
where properness is proved for convex cocompact discrete embeddings,
and Labourie~\cite{Labourie_energy}, where the above result is extended
to %% convex $\rpt$-structures.
quasi-isometric embeddings $\pi \hookrightarrow G$.
%%%% maybe say some more

\subsection{Holomorphic vector bundles and uniformization}
Let $\pi\xrightarrow{\rho}\PSLtR$ be a Fuchsian representation
corresponding to a marked hyperbolic structure $\Sigma\xrightarrow{f}M$.
A {\em spin structure\/} on $\Sigma$ determines a lifting of $\rho$
to 
\begin{equation*}
\pi\xrightarrow{\tilde\rho}\SLtR \subset \SLtC 
\end{equation*}
and hence a flat $\C^2$-bundle $(\C^2)_\rho$ over $\Sigma$. 

Choose a marked Riemann surface $X$ corresponding to a point in
Teichm\"uller space $\Teich$. 
Since locally constant maps are holomorphic for {\em any\/} complex structure
on $\Sigma$, the flat bundle $(\C^2)_\rho$ has a natural holomorphic structure;
denote the corresponding holomorphic rank two vector bundle over $X$ by 
$E_\rho \to X$. 
% Furthermore a spin structure on $\Sigma$ determines a line bundle whose
%square is the canonical bundle $K_X$; 

In trying to fit such a structure into a moduli problem over $X$,
the first problem is that this holomorphic vector bundle is {\em unstable\/}
and does not seem susceptible to Geometric Invariant Theory techniques.
Indeed, its instability intimately relates to its role in uniformization.
Namely, the developing map
\begin{equation*}
\tM \xrightarrow{\dev} \cpo 
\end{equation*}
determines a holomorphic line bundle $L\subset \Ep$. 
Since $\deg(\Ep) = 0$, and
$\dev$ is nonsingular, 
the well-known isomorphism
\begin{equation*}
T\big(\cpo\big) \cong \Hom(\gamma,\gamma^{-1}) 
\end{equation*}
where $\gamma\to\cpo$ is the tautological line bundle implies that
\begin{equation*}
L^2 \cong K_X
\end{equation*}
and $\deg(L) = g - 1 > 0$. Therefore $\Ep$ is unstable.
In fact, $\Ep$ is a nontrivial extension
\begin{equation*}
L \longrightarrow \Ep \longrightarrow \Ep/L \cong L^{-1} 
\end{equation*}
determined by the fundamental cohomology class $\varepsilon$ in 
\begin{equation*}
H^1(X,\Hom(L^{-1},L) \cong H^1(X,K) \cong \C 
\end{equation*}
defining Serre duality. (Compare Gunning~\cite{Gunning}.)

One resolves this difficulty by changing the question. 
Replace the extension class $\varepsilon$ by 
an auxiliary holomorphic object --- a Higgs field 
\begin{equation*}
\Phi \in H^0(X;K_X\otimes \End(E)) 
\end{equation*}
for the vector bundle $E := L \oplus L^{-1}$ so that the {\em Higgs pair\/}
$(E,\Phi)$ is {\em stable\/} in the appropriate sense. In our setting
the Higgs field corresponds to
the everywhere nonzero holomorphic section of the trivial holomorphic
line bundle
\begin{equation*}
\C \cong K_X\otimes \Hom(L,L^{-1})\subset K_X\otimes\End(E).
\end{equation*}
Now the only $\Phi$-invariant holomorphic subbundle of $E$ is $L^{-1}$ which
is negative, and the pair $(E,\Phi)$ is stable.

\section{Rank two Higgs bundles}
Now we impose a conformal structure on the surface to obtain extra
structure on the deformation space $\hpgg$. As before $\Sigma$ denotes
a fixed oriented smooth surface, and $X$ a Riemann surface with a fixed
marking $\Sigma \to X$.

\subsection{Harmonic metrics}

Going from $\rho$ to $(V,\Phi)$ involves finding a {\em harmonic metric,\/}
which may be regarded as a $\rho$-equivariant harmonic map
\begin{equation*}
\widetilde{M} \xrightarrow{\widetilde h} \GLnC/\mathsf{U}(n) 
\end{equation*}
into the symmetric space $\GLnC/\mathsf{U}(n)$.
The metric $h$ determines a reduction of structure group of $\Ep$ 
from $\GLnC$ to $\Un$, giving $\Ep$ a Hermitian structure. 
Let $A$ denote the unique connection on $\Ep$ which is unitary
with respect to $h$.
The harmonic metric determines the Higgs pair 
$(V,\bar{\partial}_V,\Phi)$ as follows.
\begin{itemize}
\item The Higgs field  
$\Phi$ is the holomorphic $(1,0)$-form 
$\partial h\in \Omega^1\big(\End(V)\big)$,
where the tangent space to $\GLnC/\mathsf{U}(n)$ is identified
with a subspace of $h^*\End(V)$;
\item The holomorphic structure $d_A''$ on $V$ arises from
conformal structure $\Sigma$ and the 
Hermitian connection $A$. 
\end{itemize}
The Higgs pair satisfies the {\em self-duality equations\/}
with respect to the Hermitian metric $h$:
\begin{align}
(d_A)''(\Phi) & = 0 \notag\\ 
F(A) + [\Phi,\Phi^*] & = 0 \label{eq:selfduality} 
\end{align}
Here $F(A)$ denotes the curvature of $A$, 
and $\Phi^*$ denotes the adjoint of $\Phi$ with respect to $h$.
Conversely, Hitchin and Simpson show that every stable Higgs pair
determines a Hermitian metric satisfying \eqref{eq:selfduality}.

\subsection{Higgs pairs and branched hyperbolic structures}

Choose an integer $d$ satisfying 
\begin{equation*}
0 \le d < 2 g -2
\end{equation*}
Hitchin identifies the component $\Euler^{-1}(2-2g + d)$
with Higgs pairs $(V,\Phi)$ where 
\begin{equation*}
V = L_1 \oplus L_2 
\end{equation*}
is a direct sum of line bundles $L_1$ and $L_2$
defined as follows.
Choose a square-root $K_X^{1/2}$  of the canonical bundle
$K_X$ and let $K_X^{-1/2}$ be its inverse.
Let $D\ge 0$ be an effective divisor
of degree $d$.
Define line bundles
\begin{align*}
L_1 & := K_X^{-1/2} \otimes D  \\
L_2 & := K_X^{1/2} 
\end{align*}
% Now 
% \begin{equation*}
% \deg(L_1) - \deg(L_2) = (1-g) + d - (g-1) = 2 - 2g + d = e. 
% \end{equation*}
Define a Higgs field 
\begin{equation*}
\Phi = \bmatrix 0 & s_D \\ Q & 0 \endbmatrix
\end{equation*}
where:
\begin{itemize}
\item
$s_D$ is a % meromorphic 
holomorphic 
section of the line bundle corresponding to $D$,
which determines the component of $\Phi$ in
\begin{equation*}
K_X \otimes \Hom(L_2,L_1) \cong D  
\subset \Omega^1\big(\Sigma,\End(V)\big) \big);
\end{equation*}
\item
$Q\in H^0(\Sigma,K_X^2)$ is a 
holomorphic quadratic differential with
$ \div(Q) \ge D, $
which determines the component of $\Phi$ in
\begin{equation*}
K_X \otimes \Hom(L_1,L_2) \cong K_X^2 
\subset \Omega^1\big(\Sigma,\End(V)\big) \big).
\end{equation*}
\end{itemize}
Then $(V,\Phi)$ is a stable Higgs pair.

When $Q=0$, this Higgs bundle corresponds to the uniformization
representation.
In general, when $d = 0$, the harmonic metric 
is a diffeomorphism (Schoen-Yau~\cite{SchoenYau})
% Eels-Wood~\cite{EelsWood}) and 
$Q$ is its {\em Hopf differential.\/}

The Euler class of the corresponding representation 
equals
\begin{equation*}
\deg(L_2) - \deg(L_1) = 2 - 2g + d
\end{equation*}
\begin{thm}[Hitchin~\cite{Hitchin_selfduality}]
The component $\Euler^{-1}(2 - 2g + d)$ 
identifies with a holomorphic vector bundle over the symmetric
power $\mathsf{Sym}^d(X)$. The fiber over $D\in \mathsf{Sym}^d(X)$
is the vector space
% $3g-3-d$-dimensional complex 
\begin{equation*}
\{ Q\in H^0(X,K_X^2) \mid  \div(Q) \ge D \} \cong \C^{3(g-1)-d}
\end{equation*}
\end{thm}
\noindent
The quadratic differential $Q$ corresponds to the 
{\em Hopf differential\/} of the harmonic metric $h$.
When $Q = 0$, the harmonic metric is {\em holomorphic,\/}
and defines a developing map for a branched conformal structure,
with branching defined by $D$. 

When $e = 2-2g$, then $d=0$ and the space 
$\FrickeX$ of {\em Fuchsian representations\/}
identifies with the vector space $H^0(X,K_X^2)\cong\C^{3(g-1)}$.

\subsection{Uniformization with singularities}
McOwen~\cite{McOwen} and Troyanov~\cite{Troyanov} proved a general
uniformization theorem for hyperbolic structures with conical singularities.
Specificly, let $D = (p_1) + \dots + (p_k)$ be an effective divisor, with
$p_i\in X$. Choose real numbers $\theta_i > 0$ and introduce singularities
in the conformal structure on X by replacing a coordinate chart at $p_i$
with a chart mapping to a cone with cone angle $\theta_i$.
The following uniformization theorem describes when there is a singular
hyperbolic metric in this singular conformal structure.
\begin{thm}[McOwen~\cite{McOwen}, Troyanov~\cite{Troyanov}]
If
\begin{equation*}
2 - 2 g + \sum_{i=1}^k (\theta_i  - 2\pi) > 0,
\end{equation*}
there exists a unique singular hyperbolic surface conformal to $X$ with
cone angle $\theta_i$ at $p_i$.
\end{thm}
When the $\theta_i$ are multiples of $2\pi$, then this structure is a
branched structure (and the above theorem follows from
Hitchin~\cite{Hitchin_selfduality}). 
The moduli space of such branched conformal structures 
forms a bundle $\mathfrak{S}^d$ over $\Teich$ where the fiber over a marked
Riemann surface $\Sigma\to X$ is the symmetric power $\mathsf{Sym}^d(X)$ where
\begin{equation*}
d = \frac1{2\pi}\sum_{i=1}^k (\theta_i  - 2\pi).
\end{equation*}
The resulting {\em uniformization map\/}
\begin{equation*}
\mathfrak{S}^d \xrightarrow{\mathfrak{U}} \Euler^{-1}\big(2-2g + d\big) \subset \hpgg
\end{equation*}
is {\em homotopy equivalence.\/} 
which is not surjective, by the example in \S~\ref{sec:nobranch}.

\begin{conj}
Every representation with non-discrete image lies in the
image of $\mathfrak{U}$.
\end{conj}

\section{Split $\R$-forms and Hitchin's Teichm\"uller component}

\begin{figure}[b]
\centerline{\epsfig{figure=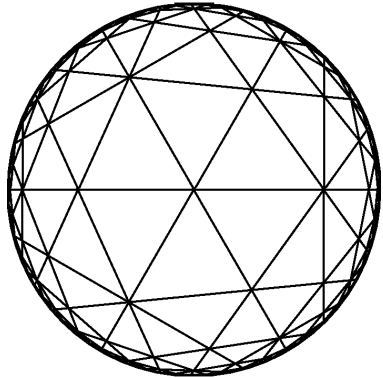,height=5cm}}
\label{fig:334}
\caption{A triangle tesselation in the hyperbolic plane, drawn in the
Beltrami-Klein projective model. Its holonomy representation is obtained
by composing a Fuchsian representation in $\SLtR$ with the irreducible
representation $\SLtR \longrightarrow \SLthR;$. }
\end{figure}
\begin{figure}[ht]
\centerline{\epsfig{figure=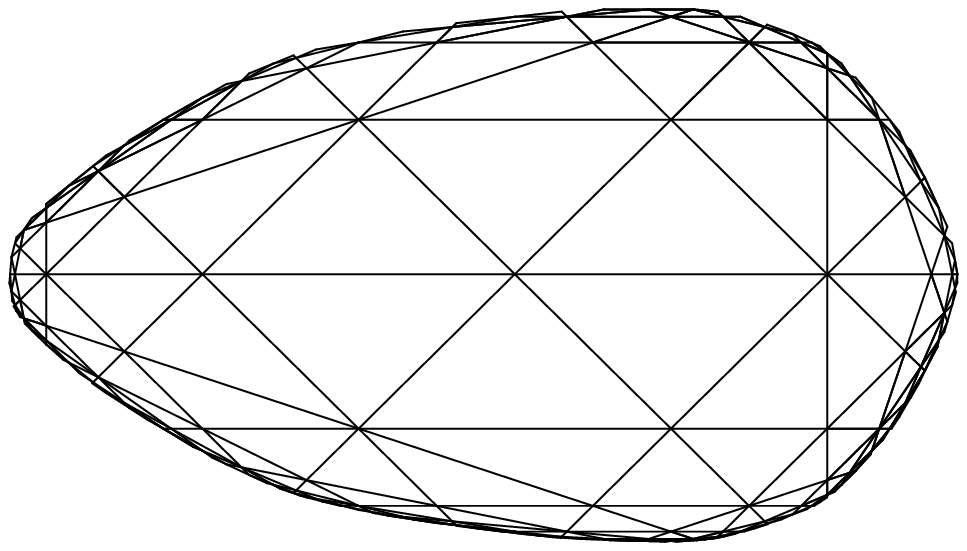,height=5cm}}
\label{fig:t334}
\caption{A deformation of a Fuchsian representation
preserving an exotic convex domain. The boundary is a strictly convex $C^1$ 
curve which is not $C^2$.}
\end{figure}

When $G$ is a split real form of a semisimple Lie group,
Hitchin~\cite{Hitchin_Lie} used Higgs bundle techniques to determine
an interesting connected component of $\hpgg$, which is {\em not detected\/}
by characteristic classes. 
A {\em Hitchin component\/} of $\hpg$ is the connected component containing
a composition
\begin{equation*}
\pi \stackrel{\rho_0}\hookrightarrow \SLtR \xrightarrow{\mathsf{K}} G
\end{equation*}
where $\rho_0$ is Fuchsian and $\mathsf{K}$ is the representation
corresponding to the {\em 3-dimensional principal subgroup\/}
discovered by Kostant~\cite{Kostant}.
When $G = \SLnR$, then Kostant's 
representation $\mathsf{K}$ is the irreducible $n$-dimensional representation
corresponding to the symmetric power $\mathsf{Sym}^{n-1}(\R^2)$.

The compositions $\mathsf{K}\circ\rho_0$ above determine a subset of $\hpgg$
which identifies with the Fricke-Teichm\"uller space, and Hitchin's main
result is that each Hitchin component is a cell of (the expected) dimension 
$\dim(G)(2g-2)$.

For example, if $G = \SLnR$, then Hitchin identifies this component with
with the $2(g-1)(n^2-1)$-cell
\begin{align*}
H^0(X;K_X^2) \oplus  
& H^0(X;K_X^3) \oplus  \dots \oplus  
H^0(X;K_X^n)  \\ & \cong
\C^{3(g-1)} \oplus \C^{5(g-1)} \oplus \dots \oplus   \C^{(2n-1)(g-1)}.
\end{align*}
When $n$ is odd, Hitchin proves there are exactly $3$ components.
The second Stiefel-Whitney characteristic class is nonzero on exactly
one component; it is zero on two components, one of which is the 
Hitchin-Teichm\"uller component.

\subsection{Convex $\rpt$-structures: $G =\SLthR$}
When $G \cong \PGLthR \cong \SLthR$, Hitchin~\cite{Hitchin_Lie}
conjectured that his component corresponded to the deformation space
$\CC$ of {\em marked convex $\rpt$-structures,\/} proved
in~\cite{Goldman_convex} to be a cell of dimension
$16(g-1)$. In~\cite{Goldman_Choi} Suhyoung Choi and the author proved
this conjecture. A {\em convex $\rpt$-manifold\/} is a quotient
$\Omega/\Gamma$ where $\Omega\subset\rpt$ is a convex domain and
$\Gamma$ a discrete group of collineations acting properly and freely
on $\Omega$.  If $\chi(M) < 0$, then necessarily $\Omega$ is {\em
properly convex\/} (contains no complete affine line), and its
boundary $\partial\Omega$ is a $C^{1+\alpha}$ strictly convex curve,
for some $0 < \alpha \le 1$. Furthermore $\alpha = 1$
if and only if $\partial\Omega$ is a conic and the $\rpt$-structure
arises from a hyperbolic structure. These facts are due to
Kuiper~\cite{Kuiper} and Benz\'ecri~\cite{Benzecri} and have recently been
extended and amplified to compact quotients of convex domains in $\rpn$
by Benoist~\cite{Benoist1,Benoist2}.

\subsection{Higgs bundles and affine spheres}
The Higgs bundle theory of Hitchin~\cite{Hitchin_Lie} identifies, for
an arbitrary Riemann surface $X$, the Hitchin component $\CC$ with the complex
vector space
\begin{equation*}
H^0(X,K_X^2) \oplus H^0(X,K_X^3)  \cong \C^{8g-8}
\end{equation*}
and the component in $H^0(X,K_X^2)$ of the Higgs field 
corresponds to the Hopf differential 
of the harmonic metric. Using the theory of {\em hyperbolic affine spheres\/}
developed by Calabi, Loewner-Nirenberg, Cheng-Yau, Gigena, Sasaki, Li,
and Wang, Labourie~\cite{Labourie_affine,Labourie_margulis} 
and Loftin~\cite{Loftin1} proved:

\begin{thm}\label{thm:LL}
The deformation space $\CC$ naturally identifies
with the holomorphic vector bundle over $\Teich$ whose fiber over a marked
Riemann surface $\Sigma\to X$ is $H^0(X,K_X^3)$..
\end{thm}
\noindent 
For every such representation, there exists a unique conformal structure so that 
\begin{equation*}
\tilde\Sigma  \xrightarrow{\tilde h} \SLthR/\mathsf{SO}(3)
\end{equation*}
is a {\em conformal map,\/} that is the component of the Higgs field
in $H^0(\Sigma,K_X^2)$  --- the Hopf differential $\Hopf(h)$ --- vanishes.
This defines the projection $\CC\to\Teich$. The zero-section corresponds
to the {\em Fuchsian\/} $\rpt$-structures, that is, the $\rpt$-structures
arising from hyperbolic structures on $\Sigma$.

It is natural to attempt to generalize this as follows. 
For any split real form $G$, and Riemann surface $X$ with $\pi_1(X) \cong \pi$,
Hitchin~\cite{Hitchin_Lie} identifies a certain direct sum of holomorphic line bundles
$\mathfrak{V}_X$ naturally associated to $X$ so that
a Hitchin component of $\Hom(\pi,G)/G$ identifies with
the complex vector space
\begin{equation*}
H^0(X,K_X^2) \oplus  H^0(X,\mathfrak{V}_X).
\end{equation*}
However, this identification depends crucially on the Riemann surface $X$ and fails
to be $\Mod$-invariant. Generalizing the Labourie-Loftin Theorem~\ref{thm:LL}, 
we conjecture that each Hitchin component of $\Hom(\pi,G)/G$ identifies
naturally with the total space of a holomorphic vector bundle $\mathfrak{E}(\Sigma)$ over
$\Teich$, whose fiber over a marked Riemann surface $X$ equals $H^0(X,\mathfrak{V}_X)$.

\subsection{Hyperconvex curves}
In 2002, Labourie~\cite{Labourie_anosov} discovered an important 
property of the Hitchin component:

\begin{thm}[Labourie]
A representation in the Hitchin component for $G = \SLnR$
is a discrete quasi-isometric embedding
\begin{equation*}
\pi \stackrel{\rho}\hookrightarrow \SL(n,\R)
\end{equation*}
with reductive image. 
\end{thm}
A crucial ingredient in his proof is the following notion.
% Also true for maximal representations.
A curve $S^1 \xrightarrow{f} \rpn$ is {\em hyperconvex\/} 
if and only if for all 
$x_1, \dots, x_n\in S^1$  distinct,
\begin{equation*}
f(x_1) +  \dots + f(x_n)  = \R^n.
\end{equation*}
\begin{thm}
[Guichard~\cite{Guichard_2005,Guichard_2006}, Labourie~\cite{Labourie_anosov}]
$\rho$ is Hitchin if and only if
$\rho$ 
preserves hyperconvex
curve.
\end{thm}
Recently Fock and Goncharov~\cite{Fock_Goncharov,Fock_Goncharov_convex} 
have studied this
component of representations, using global coordinates generalizing
Thurston and Penner's {\em shearing coordinates.\/} In these coordinates
the Poisson structure admits a particularly simple expression, leading
to a quantization. Furthermore they find a {\em positive structure\/}
which leads to an intrinsic characterization of these semi-algebraic
subsets of $\hpgg$. Their work has close and suggestive connections with
cluster algebras and $K$-theory.

\section{Hermitian symmetric spaces: Maximal representations}
We return now to the maximal representations into groups of Hermitian
type, concentrating on the unitary groups $\Upq$ and the symplectic
groups $\SpnR$.

\subsection{The unitary groups $\Upq$}
The Milnor-Wood inequality \eqref{eq:MilnorWood} may be the first
example of the {\em boundedness\/} of a cohomology class.
In a series of papers~\cite{Burger_Monod_GAFA,Burger_Iozzi_complex,
Burger_Iozzi_app,Iozzi_ern,Burger_Iozzi_Wienhard_ann,
Burger_Iozzi_Wienhard_toledo},
Burger, Monod, Iozzi and Wienhard place the local and global rigidity
in the context of the Toledo invariant being a 
{\em bounded cohomology class.\/} 
A consequence of these powerful methods for surface groups is the following,
announced in \cite{Burger_Iozzi_Wienhard_ann}:

\begin{thm}[Burger--Iozzi--Wienhard~\cite{Burger_Iozzi_Wienhard_ann}]
Let $X$ be a Hermitian symmetric space, and maximal
representation
\begin{equation*}
\pi \xrightarrow{\rho} G. 
\end{equation*}
\begin{itemize}
\item The Zariski closure $L$ of $\rho(\pi)$ is reductive;
\item The symmetric space associatied to $L$ is a Hermitian symmetric
tube domain, totally geodesicly embedded in the symmetric space of $G$;
\item $\rho$ is a discrete embedding.
\end{itemize}
Conversely, if $X$ is a tube domain, then there exists a
maximal $\rho$ with $\rho(\pi)$ Zariski-dense.
\end{thm}
For example, if $G=\Upq$, where $p\le q$, then $\rho$ is conjugate to the
normalizer $\Uppq$ of $\Upp$ in $\Upq$. As in the rank one case
(compare \S\ref{sec:Toledo}), the components of maximal representations
have strictly smaller dimension.
(In earlier work Hernandez~\cite{Hernandez} considered the
case of $\mathsf{U}(2,q)$.) 

Furthermore every maximal representation deforms into the composition
of a Fuchsian representation $\pi\xrightarrow{\rho}\SUoo$ with the
diagonal embedding
\begin{equation*}
\SUoo \subset \Uoo \stackrel{\Delta}\hookrightarrow
%\stackrel{p}
\overbrace{\Uoo \times \dots \times \Uoo}^p \subset \Upp\subset \Upq
\end{equation*}

At roughly the same time, Bradlow, Garc\'ia-Prada and
Gothen~\cite{Bradlow_GarciaPrada_Gothen} investigated the space
of Higgs bundles using infinite-dimensional Morse theory, in a similar
way to Hitchin~\cite{Hitchin_selfduality}. Their critical point analysis
also showed that maximal representations formed components of strictly
smaller dimension. They found that the number of connected components of
$\hp{\Upq}$ equals:
\begin{equation*}
2 (p+q)\, \mathsf{min}(p,q)\, (g-1) \;+\; \mathsf{gcd}(p,q). 
\end{equation*}
(For a survey of these techniques and other results, compare
\cite{Bradlow_GarciaPrada_Gothen2} as well as their recent column
\cite{Bradlow_GarciaPrada_Gothen_whatis}.)

\subsection{The symplectic groups $\SpnR$}
The case $G = \Sp(2n,\R)$ is particularly interesting, since $G$ is both
$\R$-split and of Hermitian type. Gothen~\cite{Gothen} showed there
are $3\cdot 2^{2g} + 2g -4$ components of maximal representations when
$n = 2$. For $n>2$, there are $3\cdot 2^{2g}$ components of maximal representations 
Garc{\'{\i}}a-Prada, Gothen, and Mundet~i Riera~\cite{GarciaPrada_Gothen_Mundet}).
For $n=2$, the components the nonmaximal representations are just the preimages
of the Toledo invariant, comprising $1 + 2 (2g-3) = 4g-5$ components.
Thus the total number of connected components of $\Hom\big(\pi,\Sp(4,\R)\big)$
equals 
\begin{equation*}
2 \big( 3\cdot 2^{2g} + 2 g -4 \big) + 4 g - 5 \;=\; 6 \cdot 4^g + 10 g - 13.  
\end{equation*}

The Hitchin representations are maximal and comprise $2^{2g+1}$ of these
components. They correspond to deformations of compositions
of Fuchsian representations $\pi \xrightarrow{\rho_0}\SLtR$  with the
{\em irreducible\/} representation
\begin{equation*}
\SLtR \longrightarrow 
\Aut\big(\mathsf{Sym}^{2n-1}(\R^2)\big) \hookrightarrow \Sp(2n,\R)
\end{equation*}
where $\R^{2n} \cong \mathsf{Sym}^{2n-1}(\R^2)$ with the symplectic
structure induced from $\R^2$.

Another class of maximal representations 
arises from deformations of compositions of a Fuchsian
representation $\pi\xrightarrow{\rho_0}\SLtR$ with 
the {\em diagonal embedding\/}
\begin{equation*}\label{eq:diagonal}
\SLtR \stackrel{\Delta}\hookrightarrow
\overbrace{
\SLtR\times\dots\times\SLtR
}^n
\hookrightarrow
\Sp(2n,\R). 
\end{equation*}
More generally, the diagonal embedding extends to a representation
\begin{equation*}
\SLtR \times \OO(n) \stackrel{\widetilde{\Delta}}\hookrightarrow \Sp(2n,\R)
\end{equation*}
corresponding to the 
$\SLtR \times \OO(n)$-equivariant decomposition of the symplectic vector
space
\begin{equation*}
\R^{2n} = \R^2 \otimes \R^n
\end{equation*}
as a tensor product of the symplectic vector space $\R^2$ and the
Euclidean inner product space $\R^n$.  Deformations of compositions of
Fuchsian representations into $\SLtR \times \OO(2)$ with
$\widetilde{\Delta}$ provide $2^{2g}$ more components of maximal
representations.

For $n>2$, these account for all the maximal components.
This situation is more complicated when $n=2$. 
In that case, $4g -5$ components of maximal representations
into $\Sp(4,\R)$ do not contain representations into smaller
compact extensions of embedded subgroups isomorphic to $\SLtR$.
In particular the image of every representation in such a maximal
component is {\em Zariski dense\/} in $\Sp(4,\R)$, in contrast
to the situation for $\Upq$ and $\Sp(2n,\R)$ for $n>2$. 
See Guichard-Wienhard~\cite{GuichardWienhard2} for more details.

\subsection{Geometric structures associated to Hitchin representations.}

Fuchsian representations into $\SLtR$ correspond to hyperbolic
structures on $\Sigma$, and Hitchin representations into $\SLthR$
correspond to convex $\rpt$-structures on $\Sigma$.  What geometric
structures correspond to other classes of surface group
representations?

Guichard and Wienhard~\cite{GuichardWienhard} associate to a Hitchin
representation in $\SL(4,\R)$ an $\rpthree$-structure on the {\em unit
tangent bundle\/} $T_1(\Sigma)$ of a rather special type.  The trajectories
of the geodesic flow on $T_1(\Sigma)$ (for any hyperbolic metric on $\Sigma$),
develop to projective lines.  The leaves of the weak-stable foliations
of this structure develop into convex subdomains of projective planes
in $\rpthree$. The construction of this structure uses the hyperconvex
curve in $\rpthree$. This {\em convex-foliated structure\/} is a geometric 
structure corresponding to Hitchin representations in $\SL(4,\R)$. 

For the special case of Hitchin representations into $\Sp(4,\R)$
(which are readily Hitchin representations into $\SL(4,\R)$), the
convex-foliated structures are characterized by a duality. Furthermore
the symplectic structure on $\R^4$ induces a contact structure on
$T^1(\Sigma)$ which is compatible with the convex-foliated
$\rpthree$-structure.  In addition, another geometric structure on
another circle bundle over $\Sigma$ arises naturally, related to the local
isomorphism $\Sp(4,\R)\longrightarrow \OO(3,2)$ and the identfication of
the Grassmannian of Lagrangian subspaces of the symplectic vector
space$\R^4$ with the conformal compactification of Minkowski
$(2+1)$-space (the {\em $2+1$-Einstein universe.\/} 
(Compare~\cite{BCDGM} for an exposition of this geometry.)  
The interplay between the contact $\rpthree$-geometry, 
flat conformal Lorentzian structures, the dynamics of 
geodesics on hyperbolic surfaces, and the resuting deformation
theory of  promises to be a fascinating extension of ideas rooted in the work
of Nigel Hitchin.

\end{document}